\documentclass[letterpaper, 10 pt, conference]{ieeeconf} 
\IEEEoverridecommandlockouts                          

\newcommand\scalemath[2]{\scalebox{#1}{\mbox{\ensuremath{\displaystyle #2}}}} 

\usepackage{amsthm}
\usepackage{graphics} 
\usepackage{epsfig} 
\usepackage{times} 
\usepackage{amsmath} 
\usepackage{amssymb}  

\usepackage[maxbibnames=99,giveninits=true]{biblatex}
\usepackage{bm}
\usepackage{xcolor}
\usepackage{graphicx}
\usepackage{algorithm}
\usepackage[noend]{algpseudocode}
\usepackage{booktabs}
\usepackage{siunitx}
\usepackage{hyperref}

\usepackage{caption}
\usepackage{lipsum}
\usepackage{multirow}
\usepackage{diagbox} 
\usepackage{subcaption}
\usepackage[normalem]{ulem}

\usepackage[symbol]{footmisc} 

\usepackage[T1]{fontenc} 
\newcommand{\textincon}[1]{%
{\fontfamily{zi4}\selectfont #1}}

\makeatletter 
\newcommand{\vast}{\bBigg@{4}}
\newcommand{\Vast}{\bBigg@{7}}
\makeatother

\newcommand{\R}{\mathbb{R}}
\newcommand{\E}{\mathbb{E}}
\newcommand{\Pro}{\mathbb{P}}

\newcommand{\Q}{\mathbb{Q}}
\newcommand{\W}{\mathcal{W}}
\newcommand{\VAR}{{\textrm{\large \textincon{V$@$R}}}}
\newcommand{\AVAR}{{\textrm{\large \textincon{AV$@$R}}}}
\newcommand{\tran}{\textrm{\textincon{T}}}
\newcommand{\bx}{\bm{x}}
\newcommand{\bu}{\bm{u}}

\newcommand{\by}{\bm{y}}
\newcommand{\bz}{\bm{z}}
\newcommand{\bxi}{\bm{\xi}}
\newcommand{\bw}{\bm{w}}

\newtheorem{definition}{\bf Definition}
\newtheorem{assumption}{\bf Assumption}
\newtheorem{theorem}{\bf Theorem}
\newtheorem{lemma}{\bf Lemma}

\title{\LARGE \bf
 Data-Driven Distributionally Robust Mitigation of \\
 Risk of Cascading Failures
}

\author{Guangyi Liu, Arash Amini, Vivek Pandey, and Nader Motee
\thanks{
G. Liu, V. Pandey, and N. Motee are with the Department of Mechanical Engineering and Mechanics, Lehigh University, Bethlehem, PA, 18015, USA. {\tt\small \{gliu,vkp219,motee\}@lehigh.edu}.\endgraf
A. Amini is with the Department of Aerospace Engineering, University of Texas at Austin. {\tt\small a.amini@utexas.edu}
}
}

\begin{document}

\maketitle

\thispagestyle{plain}
\pagestyle{plain}

\begin{abstract}
We introduce a novel data-driven method to mitigate the risk of cascading failures in delayed discrete-time Linear Time-Invariant (LTI) systems. Our approach involves formulating a distributionally robust finite-horizon optimal control problem, where the objective is to minimize a given performance function while satisfying a set of distributionally chances constraints on cascading failures, which accounts for the impact of a known sequence of failures that can be characterized using nested sets. The optimal control problem becomes challenging as the risk of cascading failures and input time-delay poses limitations on the set of feasible control inputs. However, by solving the convex formulation of the distributionally robust model predictive control (DRMPC) problem, the proposed approach is able to keep the system from cascading failures while maintaining the system's performance with delayed control input, which has important implications for designing and operating complex engineering systems, where cascading failures can severely affect system performance, safety, and reliability.
\end{abstract}


\section{Introduction}
Classical control theory is primarily concerned with designing controllers to stabilize a given system while optimizing its performance index in the presence of uncertainties in the system model and disturbances that affect the system. However, the design methods used in classical control theory rely on the assumption that the probability distribution of the system's parameters and disturbances is known, typically assuming a Gaussian distribution \cite{zhou1998essentials,doyle2013feedback}. In practical applications, relying on this assumption can be limiting since obtaining accurate probability distributions can be a challenging task. Distributionally robust optimization (DRO) methods provide an effective approach to handling uncertainties in the probability distribution of the system's parameters and disturbances. By minimizing the worst-case expected value of the objective function over a set of possible probability distributions, DRO methods can design robust controllers that can handle a wide range of uncertainties, including those that are difficult to quantify explicitly \cite{rahimian2019distributionally,shapiro2022conditional,kuhn2019wasserstein}. The DRO approaches are beneficial when the probability distributions are not well known or difficult to estimate accurately. Moreover, it offers a principled approach to balancing the trade-off between performance and robustness, enabling decision-makers to make informed trade-offs that account for the degree of uncertainty in the system. As a result, incorporating distributionally robust methods in control theory can enhance the reliability and robustness of control systems, particularly in scenarios with uncertainties in both system parameters and disturbances.

Model Predictive Control (MPC) \cite{camacho2013model}, including its stochastic and robust variants \cite{mesbah2016stochastic,cannon2005optimizing} involves optimizing a cost function subject to constraints on the system's inputs and states over a finite time horizon. As an extension of MPC, Distributionally Robust Model Predictive Control (DRMPC) aims to find a robust control policy for a set of possible probability distributions that describe the system's uncertainties. Recent research \cite{micheli2022data,coppens2021data,li2021distributionally, zolanvari2022data,mark2020stochastic,zhong2021data} has focused on developing DRMPC algorithms that can handle uncertainties in the system's parameters and disturbances, which provides improved performance and robustness in complex and uncertain systems.

In addition to the uncertainty of the disturbances, cascading failures and system time-delay are essential phenomenons to study since they can have widespread and severe consequences in control systems, ranging from power grids \cite{somarakis2021risk} and transportation networks \cite{Somarakis2020b,liu2021risk} to financial systems \cite{smolyak2020mitigation} and social networks \cite{yi2015modeling,chen2012time}. A cascading failure occurs when a localized disruption in one part of the system triggers a chain reaction of failures propagating through the system, leading to many more significant and widespread failures \cite{liu2022risk}. When designing a robust controller, it is critical to possess a plan when some particular undesired event has occurred in the system and may trigger cascading failures. In addition, designing the optimal control w.r.t the cascading failure reveals potential consequences of disruptions in the system and develops effective contingency plans to minimize the risk of such events. 

We aim to construct an effective DRMPC approach that accounts for the risk of cascading failure \cite{liu2021risk}, time-delay, and uncertainty in the disturbances. By converting the system dynamic into a compact form and imposing the condition on the existing failure, the first step is to sample input-state trajectories from the noisy system, which produces a set of data points that provide information about the unknown disturbances. Then, the ambiguity set is constructed using the sampled trajectories, which is a critical component of the DRMPC formulation since it determines the level of conservatism or aggressiveness in the resulting solution. The next step is to form the distributionally robust problem, which involves finding a solution that minimizes the worst-case expected value of the objective function over the ambiguity set, as well as satisfying the risk-related safety constraints. This optimization problem is typically infinite-dimensional and computationally inefficient, so an equivalent conic convex optimization problem is formulated and solved instead. This analysis allows decision-makers to understand how the solution might change in different scenarios and make informed decisions robust to uncertainty.

{\it Our Contributions:} We introduce a novel data-driven approach for mitigating the risk of {\it cascading failures} in {\it delayed} discrete-time LTI systems. The proposed method utilizes a distributionally robust model predictive control (DRMPC) framework that minimizes a piecewise affine performance function subject to risk-related safety constraints on the Average-Value-at-Risk $\AVAR$ \cite{rockafellar2002conditional} of input and state variables considering the existence of the input time-delay. Furthermore, to mitigate the cascading failure whose severity is characterized by nested sets, the proposed approach includes additional constraints that ensure that the safety-critical observable of the failure will not worsen over the finite time horizon. Our results demonstrate that the proposed approach effectively keeps the system from cascading failures while maintaining the system's performance with delayed control input.

\section{Preliminaries}

The $n-$dimensional Euclidean space is denoted by $\R^n$ and the set of natural numbers are defined by $\mathbb{N}$. We denote the vector of all ones by $\bm{1}_n = [1, \dots, 1]^\tran$. 



{\it Risk Measure:} The notion of Average (Conditional) Value-at-Risk ($\AVAR$) measure \cite{rockafellar2002conditional} indicates the severity of a random variable landing inside an undesirable set of values that characterizes the dangerous state of the system operation with a specific confidence level, i.e., Value-at-Risk ($\VAR$). In probability space $(\Omega, \mathcal{F}, \mathbb{P})$, the $\VAR$ of the random variable $Y: \Omega \rightarrow \R$ is defined as
\begin{equation*}
    \VAR_{\alpha}(Y) = \min\{z \mid F_{Y} (z) \geq 1 - \alpha\},
\end{equation*}
where the cumulative distribution function $F_Y(z) = \mathbb{P}\{Y \leq z\}$. Then, the $\AVAR$ is defined as follows.
\begin{definition}
    The Average-Value-at-Risk with the confidence level $(1-\alpha) \in (0,1)$ is the mean of the generalized $\alpha-$tail distribution:
\begin{equation*}
    \AVAR_{\alpha} (Y) = \int_{-\infty}^{\infty} z ~ d F_{Y}^{\alpha}(z),
\end{equation*}
where
\begin{equation*}
    F_{Y}^{\alpha}(z) = 
    \begin{cases}
        0, \hspace{2.35cm} \text{if } z < \VAR_{\alpha}(Y)\\
        \frac{F_{Y}(z) + \alpha - 1}{\alpha}, \hspace{1cm} \text{if } z \geq \VAR_{\alpha}(Y)
    \end{cases}.
\end{equation*}
\end{definition}
A smaller value of $\alpha$ indicates a higher level of confidence in the random variable $Y$ to stay below  $\VAR_{\alpha}(Y)$. When $Y$ has a continuous distribution function, $\AVAR$ can be obtained as the conditional expectation of $Y$ subject to $Y \geq \VAR_{\alpha}(Y)$.

{\it Systemic Event and Systemic Set:} A {\it systemic} event is a failure that will potentially lead to an overall malfunction of the network \cite{Somarakis2020b,liu2021risk}. In probability space $(\Omega, \mathcal{F}, \mathbb{P})$, the set of systemic events of random variable $y: \Omega \rightarrow \R$ is defined as $\{ \omega \in \Omega ~|~y(\omega) \in \W^*\}$. We define a collections of supersets $\{\W_{\delta}~|~\delta \in [0,\infty]\}$ of $\W^*$ that satisfy the following conditions for any sequence $\{\delta_n\}^{\infty}_{n=1}$ with property $\lim_{n \rightarrow \infty} \delta_n = \infty$ 
\begin{itemize}
    \item $\W_{\delta_2} \subset \W_{\delta_1}$ when $\delta_1 < \delta_2$
    \item $\lim_{n \rightarrow \infty} \W_{\delta_n} = \bigcap_{n=1}^{\infty} \W_{\delta_n} = \W^*$.
\end{itemize}

Suppose that the knowledge of a soft failure is given priory such that the severity of the corresponding event can be represented by $\delta = \sup \{\delta ~| ~ y \in \mathcal{W}_{\delta}\}$, the value of $\delta$ implies how dangerously the current state $y$ is close to the systemic event. Using the systemic level sets to represent the existing soft failures allows for a more systematic and rigorous risk assessment and management approach. Furthermore, by defining different safety levels, one can evaluate the trade-offs between system performance and risk based on the specific application and requirements.

\section{Problem Statement}\label{sec:problem_statement}

Let us consider a delayed, discrete-time, and linear time-invariant (LTI) system, which is governed by the following dynamics
\begin{equation}    \label{eq:dyn}
    \bx_{t+1} = A \bx_{t} + B \bu_{t-\tau} + \bw_t,
\end{equation}
where $\bx_t = [x^{(1)}_t,...,x^{(n)}_t] \in \R^n$ denotes the state of the system, $\bu_t = [u^{(1)}_t,...,u^{(m)}_t] \in \R^m$ denotes the control input, and $\bw_t \in \R^n$ represents the i.i.d. additive disturbance which follows some {\it unknown} probability distribution $\Pro_{w}$. The state and the input matrix are given by $A \in \R^{n\times n}$ and $B \in \R^{n\times m}$, respectively. The constant $\tau \in \mathbb{N}$ represents the input time-delay, which is commonly considered in the applications of networked control systems \cite{Somarakis2020b,liu2021risk}.

We consider the model predictive control problem (MPC) of the system \eqref{eq:dyn} for $T$ future steps with the initial states of the system denoted by $\bx_0 = [x^{(1)}_0,...,x^{(n)}_0] \in \R^n$, and the initial inputs are given by $\bu_{init} = [\bu_{-\tau},...,\bu_0]^\tran$ due to the existence of the time-delay. For the notation simplicity, we introduce the $T$-steps compact form \cite{micheli2022data} of \eqref{eq:dyn}, which contains the information of future $T$-stages of the system
\begin{equation} \label{eq:compact-dyn}
    \by = \bar{A} \bz + \bxi,
\end{equation}
where the collective states are shown by $\by = [\bx_1,\bx_2,$ $\hdots, \bx_T]^\tran \in \R^{nT}$. Vector $\bz$ represents the comprehensive information of both initial conditions and control inputs, i,e, $\bz = [\bx_0^\tran,\bu_{init}^\tran, \bu_e^\tran]^\tran \in \R^{n+mT}$, where the effective control input is given by
$
    \bu_e := [\bu_{1},...,\bu_{T-1-\tau}]^\tran.
$
The collective disturbances are given by $\bxi = \bar{B} \bw$, where $\bw = [\bw_0, \bw_1, \hdots, \bw_{T-1}]^\tran$. Matrices $\bar{A}$ and $\bar{B}$ are defined as follows:
\begin{align*}
    \bar{A} = \scalemath{0.68}{\begin{bmatrix} 
    A & B & 0_{n\times m} & \hdots & 0_{n\times m}\\ 
    A^2 & AB & B & \ddots & \vdots \\ 
    \vdots & \vdots & \ddots & \ddots & 0_{n\times m}\\ 
    A^T & A^{T-1}B & \hdots & AB & B
    \end{bmatrix}}, 
    \bar{B} = \scalemath{0.68}{\begin{bmatrix} 
    I_n & 0_{n\times n} & \hdots & 0_{n\times n}\\ 
    A & I_n & \ddots & \vdots \\ 
    \vdots & \ddots & \ddots & 0_{n\times n}\\ 
    A^{T-1} & \hdots & A & I_n
    \end{bmatrix}}.
\end{align*}

Let us also introduce an observable $h_t(\by,\bz)$, which represents the safety-critical states of the system at each time instance \footnote{$h_0(\by,\bz):= h_0(\bx_0)$ since it only depends on the initial condition.}, e.g., the inter-vehicle distance of a platoon of cars \cite{Somarakis2020b}. We consider the scenario when the functioning system is observed with a soft failure at $t =0$, i.e., $h_0(\bx_0) \in \W_\delta$, which is dangerously close to the systemic event with level $\delta$. To keep the system running effectively and safely, one has to design a sequence of control inputs that maintain the performance of the system while mitigating the risk of cascading soft failures. It is also worth noticing that the probability distribution of $\by$ is a function of both initial condition and control inputs, which is shown by $\Pro_{|\delta,\bz}$, where $\delta$ denotes the severity of the existing soft failure on $h_0(\bx_0)$. Then, the corresponding MPC problem is constructed as follows
\begin{equation} \label{eq:ori_mpc}
    \begin{aligned}
        &J := &\text{minimize}_{\bz}  &~\E^{\by \sim \Pro_{|\delta,\bz}} \left[c(\by,\bz)\right]\\
        &&\text{subject to}   &~\AVAR^{\by \sim \Pro_{|\delta,\bz}}_{\alpha} \left( g(\by,\bz) \right) \leq 0\\
        && &~\AVAR^{\by \sim \Pro_{|\delta,\bz}}_{\alpha} \left( h_t(\by,\bz) \right) \notin \mathcal{W}_{\delta}\\
        &&& ~ \forall t = 1, ..., T,
    \end{aligned}
\end{equation}
where $c(\by,\bz)$ defines the collective system performance measure or $T$ steps, and $g(\by,\bz)$ is the function that measures the overall safety of the system, and $h_1(\by,\bz), ..., h_T(\by,\bz)$ are the future states of the safety-critical observable $h_0(\bx_0)$. The $\AVAR$ constraints on $h_1(\by,\bz), ..., h_T(\by,\bz)$ are ensuring that the existing soft failure on the safety-critical observable $h_0(\bx_0)$ will not worsen in the future, e.g., the distance between a pair of vehicles close to collision, and the $\AVAR$ constraint on $g(\by,\bz)$ is introduced to keep the system operating away from the potential cascading failures induced by $h_0(\bx_0) \in \W_\delta$, e.g., the distance between another pair of vehicles within the same platoon.

Our {\it objective} is to solve a distribution robust optimization problem to obtain a feasible solution for problem \eqref{eq:ori_mpc} that can handle worst-case scenarios for both cost and safety constraints since information about the disturbance $\bxi$ is lacking. Additionally, the solution to problem \eqref{eq:ori_mpc} also mitigates the impact of input time-delay $\tau$ and potential cascading failures triggered by the existing soft failure $h_0(\bx_0) \in \W_\delta$.

\section{Uncertainty of the Disturbance} \label{sec:uncertainty}

\subsection{Disturbance Distribution Approximation}

To account for the unknown additive input noise $\bw_t$, it is necessary to gather sufficient information by collecting multiple sample trajectories of the system \eqref{eq:dyn}. In this study, we conduct $N$ test trajectories, each comprising a sequence of input $\tilde{\bz}$ \footnote{The input sequence includes the initial condition $\bx_0$ and initial input $\bu_{init}$.} with length $T$, $\delta = 0$, and $\tau = 0$. The input sequences $\tilde{\bz}^i$ are randomly generated, and the compact system produces corresponding outputs $\tilde{\by}^i$ for $i = 1,...,N$. The noise residual can then be calculated as a function of the input and output measurements as follows \cite{micheli2022data,verweij2003sample}
\begin{equation} \label{eq:noise_dyn}
    \tilde{\xi}^i := \tilde{\by}^i - \bar{A} \tilde{\bz}^i. 
\end{equation}
Then, for given initial condition $\bx_0$ with existing soft failure $h_0(\bx_0) \in \W_\delta$, initial input $\bu_{init}$, and effective input $\bu_e$, the $T$-step state approximation $\bar{\by}^i$ of the system \eqref{eq:compact-dyn} can be computed as 
\begin{align*}
    \bar{\by}^i := \bar{A} \bz + \tilde{\bxi}^i, \text{ with } \bz = [\bx_0^\tran,\bu_{init}^\tran, \bu_e^\tran]^\tran,
\end{align*}
and its empirical distribution is given by
\begin{align}
    \tilde{\Pro}_{|\delta,\bz} := \frac{1}{N}\sum_{i =1}^{N} \bm{\delta}({\bar{\by}^i}),
\end{align}
where $\bm{\delta}(\by)$ is the Dirac distribution at $\by$.

\subsection{Wasserstein Distance and the Ambiguity Set}

In the sampling average approximation process for obtaining the distribution of the disturbance $\bxi$ \cite{verweij2003sample}, due to its inherent nature of uncertainty, it is nearly impossible to obtain the exact distribution of $\bxi$. Consequently, there exists a discrepancy between the true distribution and the approximated distribution of the $T$-step state $\by$ for some given $\bz$, denoted as $\Pro_{|\delta,\bz}$ and $\tilde{\Pro}_{|\delta,\bz}$, respectively. This difference raises questions regarding how to solve the MPC problem with $\Tilde{\Pro}_{|\delta,\bz}$ while ensuring the performance and safety of the system. The Wasserstein distance can be used to quantify this discrepancy, as it is a commonly employed metric for measuring the difference between two potentially distinct probability distributions \cite{mohajerin2018data,kantorovich1958space}, which is defined on the space $\mathcal{M}(\Xi)$ of all probability distributions $\Q$ supported on $\Xi$ such that $\E_{\Q} [|\xi|] = \int_{\Xi} |\xi| \Q (\text{d} \xi) < \infty$.

\begin{definition}
    (Wasserstein distance). For any $p \in [1, \infty)$, the type-$p$ Wasserstein distance between two probability distributions $\Pro$ and $\Q$ supported on $\R^n$ is defined as 
    \begin{equation}
        W_p\left(\Pro, \Q\right) := \left( \inf_{\pi \in \Pi(\Pro, \Q)} \int_{\R^n \times \R^n} \|\xi -\xi'\|^p \pi(\text{d} \xi, \text{d} \xi')  \right)^{\frac{1}{p}},
    \end{equation}
    where $\|\cdot\|$ is a norm on $\R^n$, while $\Pi(\Pro, \Q)$ denotes the set of all joint probability distributions of $\xi \in \R^n$ and $\xi' \in \R^n$ with marginals $\Pro$ and $\Q$, respectively.
\end{definition}

The Wasserstein distance is a metric in $\R^n$ that satisfies several properties, including non-negativity, symmetry, sub-additivity, and the fact that it vanishes if and only if $\Pro = \Q$ \cite{villani2009optimal}. Moreover, if both $\Pro$ and $\Q$ have finite $p$'th order moments, then $W_p\left(\Pro, \Q\right)$ is finite. As a result, the Wasserstein distance can be used to construct an ambiguity set of probability distributions.

\begin{definition}(Ambiguity set)
    For any $\Pro$ and $\Q$ that is supported on $\R^{n}$ such that $\E^{\Pro}[\|\by\|] < \infty$ and $\E^{\Q}[\|\by\|] < \infty$, the ambiguity set centered at $\Pro$ with radius $r$ is defined as
    \begin{equation}
    \mathfrak{M} (\Pro):= \left\{\mathbb{Q} ~|~ W_p(\Pro, \Q) \leq r \right\},
    \end{equation}
    where $\Pro$ is the target probability distribution, and $r > 0$ defines the radius of the ambiguity set (Wasserstein ball).
\end{definition}

\begin{figure}[t]
    \centering
    \includegraphics[width = \linewidth,height = 3cm]{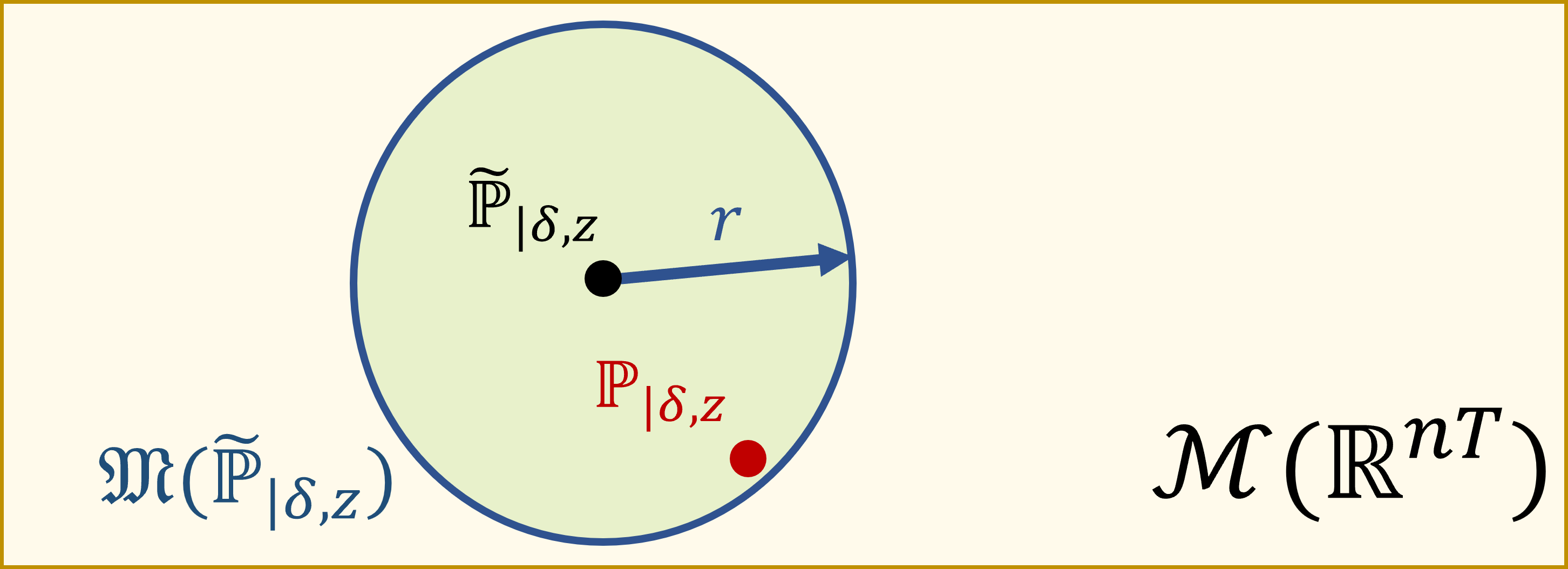}
    \caption{The ambiguity set $\mathfrak{M}(\tilde{\Pro}_{|\delta,\bz})$ with radius $r$ that captures the true distribution $\Pro_{|\delta,\bz}$.}
    \label{fig:ambiguity_set}
\end{figure}

This study aims to create an ambiguity set that encapsulates the true distribution of $\by$, denoted by $\Pro_{|\delta,\bz}$, which is constructed using a Wasserstein ball centered at $\tilde{\Pro}_{|\delta,\bz}$, as illustrated in Fig. \ref{fig:ambiguity_set}. By doing so, we aim to transform the original MPC problem \eqref{eq:ori_mpc} into a distributionally robust optimization problem.

\subsection{Quantifying the Radius of the Ambiguity Set}

Due to the stochastic nature of the problem, it is impossible to guarantee that the true distribution $\Pro_{|\delta,\bz}$ will be captured by the ambiguity set $\mathfrak{M}(\tilde{\Pro}_{|\delta,\bz})$ with radius $r$. To address this uncertainty, it is reasonable to introduce a confidence level $\varepsilon \in (0,1)$ that reflects the likelihood of the true distribution being captured by the ambiguity set with type-$1$ Wasserstein distance, such that
\begin{align} \label{eq:dist_P_tP}
    \Pro\{\Pro_{|\delta,\bz} \in \mathfrak{M}(\tilde{\Pro}_{|\delta,\bz})\}=\Pro\{W_1(\Pro_{|\delta,\bz}, \tilde{\Pro}_{|\delta,\bz}) \leq r\} \geq \varepsilon.
\end{align}

The following assumption is introduced to quantify the radius of the ambiguity set using the Wasserstein distance with a finite sample size $N$.

\begin{assumption} \label{asp:light-tail}
    For some constants $\beta > 1$, 
    \begin{align*}
        H:= \E^{\Pro}\left[e^{\|\bm{\xi}\|^\beta}\right] < \infty.
    \end{align*}
\end{assumption}

The above assumption essentially necessitates the tail of the distribution $\Pro$ to exhibit exponential decay. The modern measure concentration result \cite{fournier2015rate,mohajerin2018data} provides the radius of the ambiguity set for a given confidence level $\varepsilon$ and number of collected samples $N$, such that \eqref{eq:dist_P_tP} holds, which is presented below.

\begin{lemma} \label{lem:1}
    If Assumption \ref{asp:light-tail} holds, for confidence level $\varepsilon \in (0,1)$, the condition \eqref{eq:dist_P_tP} will be satisfied with the radius of the ambiguity set $r$, such that
    \begin{equation} \label{eq:radius}
        r := \begin{cases}
            \left(\frac{\log\left(\frac{c_1}{1-\varepsilon}\right)}{c_2 N}\right)^{\frac{1}{nT}} &\text{ if } N \geq \frac{\log\left(\frac{c_1}{1-\varepsilon}\right)}{c_2 N}\\
            \left(\frac{\log\left(\frac{c_1}{1-\varepsilon}\right)}{c_2 N}\right)^{\frac{1}{\beta}} &\text{ if } N < \frac{\log\left(\frac{c_1}{1-\varepsilon}\right)}{c_2 N}
        \end{cases},
    \end{equation}
    where $c_1$ and $c_2$ are positive constants that only depend on $\beta, H$ and $nT$.
\end{lemma}

The result above offers guidance on selecting the appropriate radius for constructing the ambiguity set, ensuring that the true distribution $\Pro_{|\delta,\bz}$ is encompassed within the ambiguity set $\mathfrak{M}(\tilde{\Pro}_{|\delta,\bz})$ with a confidence level of $\varepsilon$.

\section{The Distributionally Robust Formulation of the MPC Problem}

The distributionally robust formulation of problem \eqref{eq:ori_mpc} necessitates the consideration of worst-case scenarios for both cost and safety constraints across all distributions within the ambiguity set $\mathfrak{M}(\tilde{\Pro}_{|\delta,\bz})$. The optimal solution of the distributionally robust problem, denoted by $\bz^*$, must fulfill the following condition for all $t = 1,...,T$:
\begin{equation} \label{eq:mpc_cond}
    \scalebox{1}{$
    \begin{aligned}
            &\E^{\by \sim \Pro_{|\delta,\bz^*}} \left[ c(\by,\bz^*)\right] \leq \sup_{\Q \in \mathfrak{M}(\tilde{\Pro}_{|\delta,\bz^*})} \E^{\by \sim \Q} \left[ c(\by,\bz^*)\right] \\
            &\AVAR^{\Pro_{|\delta,\bz^*}}_{\alpha} (g(\by,\bz^*)) \leq \sup_{\Q \in \mathfrak{M}(\tilde{\Pro}_{|\delta,\bz^*})}\AVAR^{\Q}_{\alpha} (g(\by,\bz^*)) \\
            &\AVAR^{\Pro_{|\delta,\bz^*}}_{\alpha} (h_t(\by,\bz^*)) \leq \sup_{\Q \in \mathfrak{M}(\tilde{\Pro}_{|\delta,\bz^*})}\AVAR^{\Q}_{\alpha} (h_t(\by,\bz^*))
    \end{aligned}$}.
\end{equation}

Then, one should seek a distributionally robust formulation of the original MPC problem \eqref{eq:ori_mpc}, which aims to develop a solution to withstand the scenario by considering the worst-case expectation and worst-case $\AVAR$ safety constraints within an ambiguity set. Hence, the distributionally robust model predictive control (DRMPC) problem can be written as
\begin{equation} \label{eq:dr_mpc}
    \begin{aligned}
        &\tilde{J} := &\text{minimize}_{\bz}  &\sup_{\Q \in \mathfrak{M}(\tilde{\Pro}_{|\delta,\bz})} \E^{\by \sim \Q} \left[ c(\by,\bz)\right]\\
        &&\text{subject to}   &\sup_{\Q \in \mathfrak{M}(\tilde{\Pro}_{|\delta,\bz})} \AVAR^{\by \sim \Q}_{\alpha} \left( g(\by,\bz) \right) \leq 0\\
        && &\sup_{\Q \in \mathfrak{M}(\tilde{\Pro}_{|\delta,\bz})} \AVAR^{\by \sim \Q}_{\alpha} \left( h_t(\by,\bz) \right) \notin \W_\delta  \\
        &&& \hspace{1.6cm} \forall t = 1,...,T.
    \end{aligned}
\end{equation}

Solving the above optimization problem can provide a potentially feasible solution to the original problem \eqref{eq:ori_mpc}, satisfying the condition \eqref{eq:mpc_cond}. As discussed in \S \ref{sec:uncertainty}, the decision-dependent ambiguity set $\mathfrak{M}(\tilde{\Pro}_{|\delta,\bz})$ may not fully capture the true distribution $\Pro_{|\delta,\bz}$, with a probability of $1-\varepsilon$, which implies that the following conditions equivalent to \eqref{eq:mpc_cond} should be introduced for all $t = 1,...,T$: 
\begin{equation} \label{eq:dr_prob_condition}
    \scalebox{0.9}{$
    \begin{aligned}
            &\Pro\left\{\E^{\by \sim \Pro_{|\delta,\bz^*}} \left[ c(\by,\bz^*)\right] \leq \sup_{\Q \in \mathfrak{M}(\tilde{\Pro}_{|\delta,\bz^*})} \E^{\by \sim \Q} \left[ c(\by,\bz^*)\right] \right\} \geq \varepsilon\\
            &\Pro\left\{\AVAR^{\Pro_{|\delta,\bz^*}}_{\alpha} (g(\by,\bz^*)) \leq \sup_{\Q \in \mathfrak{M}(\tilde{\Pro}_{|\delta,\bz^*})}\AVAR^{\Q}_{\alpha} (g(\by,\bz^*)) \right\} \geq \varepsilon\\
            &\Pro\Bigg\{\AVAR^{\Pro_{|\delta,\bz^*}}_{\alpha} (h_t(\by,\bz^*)) \leq \\
            &\hspace{3.5cm}\sup_{\Q \in \mathfrak{M}(\tilde{\Pro}_{|\delta,\bz^*})}\AVAR^{\Q}_{\alpha} (h_t(\by,\bz^*)) \Bigg\} \geq \varepsilon.
    \end{aligned}$}
\end{equation}

The abovementioned condition suggests that applying the optimal solution $\bz^*$ to the original MPC problem typically results in a conservative cost estimate, and it is likely feasible for the original problem \eqref{eq:ori_mpc} with a probability of $\varepsilon$. The conservativeness of the solution $\bz^*$ in the distributionally robust problem for the worst-case scenario is profoundly affected by the radius $r$ of the ambiguity set $\mathfrak{M}(\tilde{\Pro}_{|\delta,\bz})$. A larger radius in the ambiguity set typically implies that more potential probability distributions will be considered, making the solution more conservative. Moreover, a larger value of $\varepsilon$ also signifies a broader range of possible probability distributions, leading to a more conservative solution for the distributionally robust optimization problem.

\section{Solution to the Distributionally Robust MPC Problem}

Considering the fact that taking the supremum over probability distributions within the ambiguity set $\mathfrak{M}(\tilde{\Pro}_{|\delta,\bz})$ leads to an infinite-dimensional problem, which is both computationally inefficient and intractable. To overcome this issue, we introduce a tractable and equivalent formulation of the DRMPC problem by using the similar technique introduced in \cite{mohajerin2018data}. In order to convert DRMPC problem \eqref{eq:dr_mpc} into a tractable convex optimization problem, the following assumptions \cite{mohajerin2018data} on $c(\cdot),g(\cdot),h_1(\cdot), ...,h_T(\cdot)$ are introduced.

\begin{assumption} \label{asp:affine}
    Functions $c(\by,\bz)$, $g(\by,\bz)$, and real-valued function $h_1(\by,\bz), ..., $ $ h_T(\by,\bz)$ are piecewise affine functions, which can be expressed as \footnote{The values of $N_{i_c}, a_{i_c}, b_{i_c}$, and $c_{i_c}$ can be uniquely defined for a given piecewise affine function $c(\cdot)$. This can also be extended to $g(\cdot), h_1(\cdot), ..., h_T(\cdot)$}
    \begin{align*}
        c(\by,\bz) &= \textup{max}_{i_c \leq N_{i_c}} ~ a_{i_c} \by + b_{i_c} \bz + c_{i_c}\\
        g(\by,\bz) &= \textup{max}_{{i_g} \leq N_{i_g}} ~ d_{i_g} \by + e_{i_g} \bz + f_{i_g}\\
        h_t(\by,\bz) &= \textup{max}_{{i_h} \leq N_{i_h}^{(t)}} ~ \upsilon^{(t)}_{i_h} \by + \phi^{(t)}_{i_h} \bz + \chi^{(t)}_{i_h},
    \end{align*}
    for all $\by,\bz$, and $t = 1,...,T$.
\end{assumption}

For the exposition of the next result, let us also introduce the following notations
\begin{equation}
    \begin{aligned} \label{eq:lambda_theta}
        \lambda &= \textup{max}_{{i_c} \leq N_{i_c}} \|a_{i_c}\|,\\
        \theta &= \textup{max}_{{i_g} \leq N_{i_g}} \|d_{i_g}\|,\\
        \eta^{(t)} &= \textup{max}_{{i_g} \leq N_{i_g}^{(t)}} \|\upsilon_{i_g}^{(t)}\|.
    \end{aligned}
\end{equation}

\begin{theorem} \label{thm:2}
    If Assumption \ref{asp:affine} holds, the optimal solution of the DRMPC problem \eqref{eq:dr_mpc} can be obtained by solving the following convex optimization problem 
    \begin{equation} \label{eq:convex_conic}
        \begin{aligned}
            \tilde{J} ~:= \inf_{\bz,\bm{\nu}} ~ &\lambda ~ r + \frac{1}{N}\sum_{i = 1}^N \nu_i\\
            \textup{subject to}~ & a_{i_c}(\bar{A} \bz + \tilde{\bxi}^i) + b_{i_c} \bz + c_{i_c} \leq \nu_i \\
            & \forall ~ i = 1,...,N \textup{ and } \forall ~ {i_c} = 1,...,N_{i_c}\\
            & \bz \in \mathcal{Z}_{g} \bigcap \mathcal{Z}_{h_1} \bigcap \cdots \bigcap \mathcal{Z}_{h_T} 
        \end{aligned}
    \end{equation}
    where the set $\mathcal{Z}_{g}$ is defined for all $\bz \in \R^{n+mT}$ and some $s^{(0)},q_i^{(0)}$ such that following conditions hold 
    \begin{align*}
        \theta \, r + \frac{1}{N}\sum_{i = 1}^N q_i^{(0)} &\leq s^{(0)} \,\varepsilon\\
        \max\left\{d_{i_g}(\bar{A} \bz + \tilde{\bxi}^i) + e_{i_g} + f_{i_g} + s^{(0)},0\right\} &\leq q_i^{(0)},
    \end{align*}
    with $\forall i = 1,...,N$, and ${i_g} = 1,...,N_{i_g}$. Sets $\mathcal{Z}_{h_t}$ is defined for all $\bz \in \R^{n+mT}$ and some $s^{(t)},q_i^{(t)}$ such that following conditions hold
    \begin{align*}
        \eta^{(t)} \, r + \frac{1}{N}\sum_{i = 1}^N q_i^{(t)} &\leq s^{(t)} \,\varepsilon\\
        \max\left\{\upsilon^{(t)}_{i_h}(\bar{A} \bz + \tilde{\bxi}^i) + \phi^{(t)}_{i_h} + \chi^{(t)}_{i_h} + s^{(t)},0\right\} &\leq q_i^{(t)},
    \end{align*}
    with $\forall i = 1,...,N$, and ${i_h} = 1,...,N_{i_h}^{(t)}$.
    The value of $r$, $\lambda$, $\theta$, and $\eta^{(t)}$ are computed using \eqref{eq:radius} and \eqref{eq:lambda_theta}, accordingly.
\end{theorem}

Based on the above findings, we convert the previously intractable infinite-dimensional problem \eqref{eq:dr_mpc} into a computationally efficient conic convex optimization problem. By solving this problem numerically, one can obtain the optimal control input for the distributionally robust MPC problem.

\section{Case Studies}
\begin{figure}[t]
    \centering
    \includegraphics[width = \linewidth]{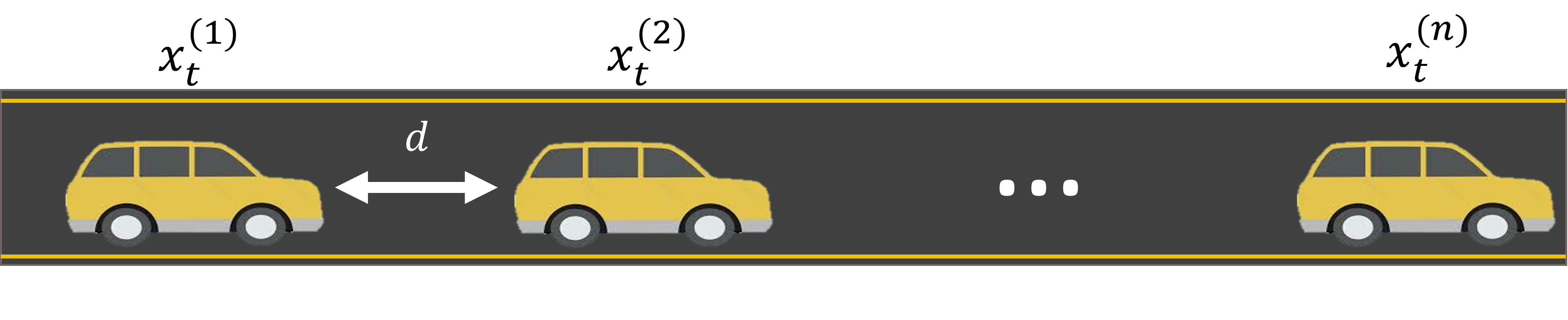}
    \caption{A team of self-driving vehicles aims to form a platoon.}
    \label{fig:platoon}
\end{figure}

\begin{figure}[t]
    \centering
    \includegraphics[width = \linewidth]{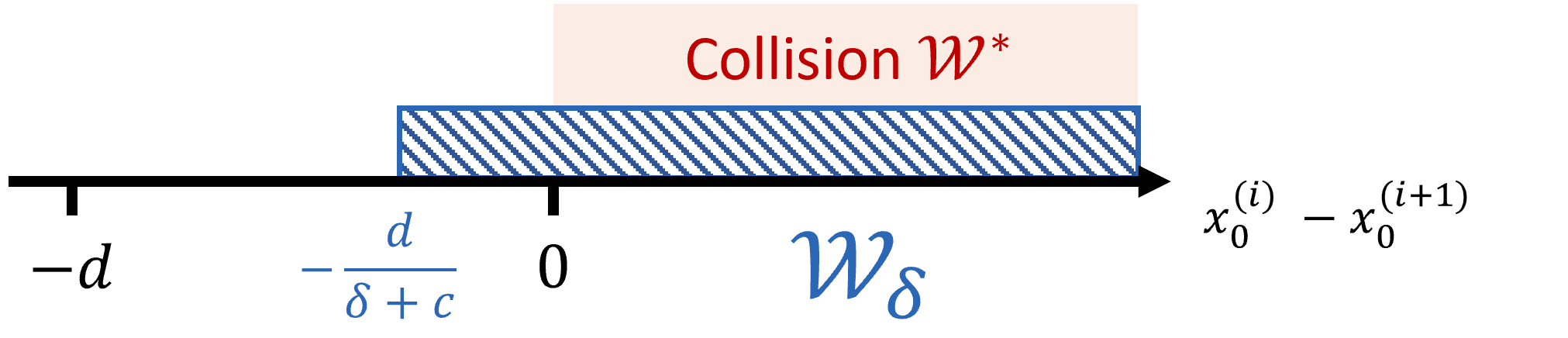}
    \caption{The above figure depicts the system event, i.e., collision, and its corresponding superset $\W_\delta$ based on the inter-vehicle distance $x^{(i)}_0 - x^{(i+1)}_0$.}
    \label{fig:W_delta}
\end{figure}

In the case study, we consider a team of $n = 6$ autonomous driving vehicles aims to form a platoon with their governing dynamics shown as follows
\begin{equation} \label{eq:case_dyn}
    \bx_{t+1} = \bx_{t} + \bu_{t-\tau} + \bw_t,
\end{equation}
where $\bx_t = [x_t^{(1)}, ...,x_t^{(n)}]^\tran \in \R^n$ denotes the position of vehicles and $\bu_{t-\tau} = [u_{t-\tau}^{(1)},...,u_{t-\tau}^{(n)}]^\tran \in \R^n$ defines their delayed control input, accordingly. The additive input noise is given by $\bw_t = [w_t^{(1)},...,w_t^{(n)}]^\tran \in \R^n$, each follows the normal distribution $\mathcal{N}(0,0.05)$. However, the distribution of the additive disturbance is considered unknown when solving the DRMPC problem. The initial condition $\bx_0$ are generated randomly and the initial input $\bu_0 = [0, ..., 0]^\tran$. The systemic event of the platoon is considered as the inter-vehicle collision, and its corresponding level set is defined as $\W_\delta = ( -\frac{d}{\delta+c}, \infty)$, with $d = 1$ and $c = 1.25$, as depicted in Fig. \ref{fig:W_delta}.
 
The {\it objective} of the system is to form a platoon such that each consecutive pair of vehicles aims to maintain a target inter-vehicle distance $d$ \cite{Somarakis2020b} for time $t \in [1,T]$ with $T=5$, as is depicted in Fig. \ref{fig:platoon}. Meanwhile, the system also aims to avoid the systemic event, i.e., the inter-vehicle collision, which will occur at the $j$'th pair of vehicles if $h_t = x^{(j)}_t - x^{(j+1)}_t> 0$. Hence, the safety-critical observable, i.e., the inter-vehicle distance at the $j$'th pair, is given by
\begin{align*}
    h_t(\by, \bz):= -\bm{\mathcal{E}}_t \, C_j \, \by,
\end{align*}
where $\{\bm{\mathcal{E}}_1, \dots, \bm{\mathcal{E}}_T\}$ denotes the set of standard Euclidean basis for $\mathbb{R}^{T}$. The matrix $C_j \in \R^{T,nT}$ is defined as
\begin{equation*}
    C_j := [\tilde{\bm{e}}_{j} \,| \, \tilde{\bm{e}}_{j+T} \,|\,... \,|\, \tilde{\bm{e}}_{j+(n-1)T}]^\tran,
\end{equation*}
where $\tilde{\bm{e}}_{j} = \bm{e}_{j+1} - \bm{e}_{j} \in \R^{nT}$, and $\{\bm{e}_1, \dots, \bm{e}_{nT}\}$ denotes the set of standard Euclidean basis for $\mathbb{R}^{nT}$. Then, the cost function is considered as 
\begin{equation}
    c(\by,\bz) = \sum_{i = 1}^{n-1}\| C_i\by - \bm{1}_{T}\|_1 + 0.04 \|\bu_e\|_1,
\end{equation}
which measures the deviation of the inter-vehicle distance w.r.t the target distance for all pairs of vehicles and all future $T$-steps and the input cost. The safety constraint functions are defined by measuring the risk of cascading collisions for all vehicles, such that
\begin{align*}
    g(\by,\bz) = -\sum_{i = 1}^{n-1}1_{T}^{\textrm{\textincon{T}}}C_i\by + \gamma_1,
\end{align*}
where $\gamma_1$ is a design parameter such that the corresponding $\AVAR$ levels satisfies the safety requirement of cascading failures. A total of $N = 50$ sample trajectories are collected from the system to approximate the distribution of the unknown additive disturbance.

The corresponding DRMPC problem is given by 
\begin{equation} \label{eq:case_mpc}
    \scalebox{0.85}{$
    \begin{aligned}
        &\tilde{J} := &\text{minimize}_{\bz} &\sup_{\Q \in \mathfrak{M}(\tilde{\Pro}_{|\delta,\bz})} \E^{\by \sim \Q} \left[ \sum_{i = 1}^{n-1}\| C_i\by - \bm{1}_{T}\|_1 + 0.04 \|\bu_e\|_1\right]\\
        &&\text{subject to}   &\sup_{\Q \in \mathfrak{M}(\tilde{\Pro}_{|\delta,\bz})} \AVAR^{\by \sim \Q}_{\alpha} \left( -\sum_{i = 1}^{n-1}1_{T}^{\textrm{\textincon{T}}}C_i\by + \gamma_1 \right) \leq 0\\
        && &\sup_{\Q \in \mathfrak{M}(\tilde{\Pro}_{|\delta,\bz})} \AVAR^{\by \sim \Q}_{\alpha} \left( h_t(\by,\bz) + \gamma_2 \right) \leq 0\\
        &&& \hspace{1.6cm} \forall t = 1,...,T,
    \end{aligned}$}
\end{equation}
where $j=1$, $\gamma_1$ is obtained from the safety design requirements, and $\gamma_2$ is computed using \eqref{eq:gamma}. The existing soft failure is measured with $h_0(\bx_0) \in \W_\delta$, indicating the $j$'th pair of vehicle is observed close to collision at $t = 0$. We solve the equivalent convex conic optimization problem to obtain the optimal control inputs $\bz^*$ and then apply them to the actual system \eqref{eq:case_dyn}.

\begin{figure}[t]
    \begin{subfigure}[t]{.32\linewidth}
        \centering
    	\includegraphics[width=\linewidth]{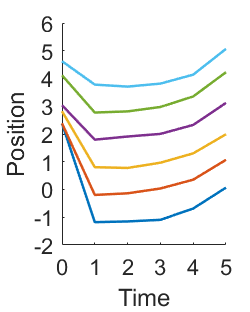}
    	\caption{$\tau = 0$.}
    \end{subfigure}
    \hfill
    \begin{subfigure}[t]{.32\linewidth}
        \centering
    	\includegraphics[width=\linewidth]{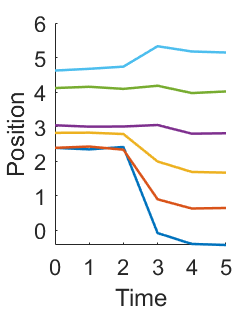}
    	\caption{$\tau = 2$.}
    \end{subfigure}
    \hfill
    \begin{subfigure}[t]{.32\linewidth}
        \centering
    	\includegraphics[width=\linewidth]{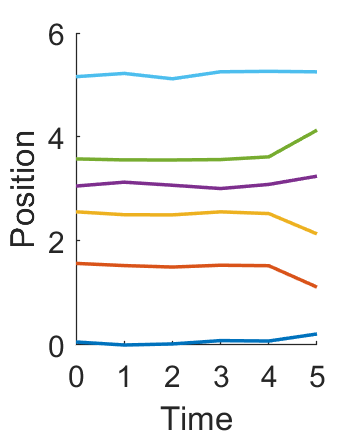}
    	\caption{$\tau = 4$.}
    \end{subfigure}
    \caption{Real time position of the platoon with random initial conditions when using the corresponding optimal control input $\bz^*$.}
    \label{fig:position_tau}
\end{figure}

\begin{figure}[t]
    \centering
    \includegraphics[width=\linewidth]{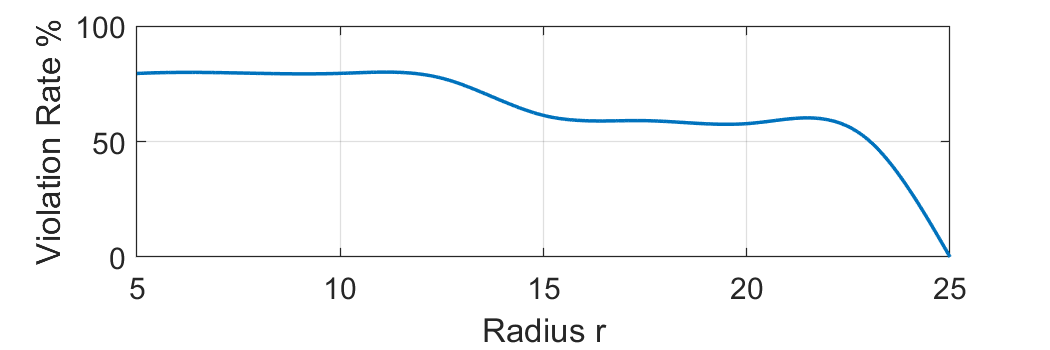}  
    \caption{The violation rate of safety constraints with various radius of the ambiguity set.}
    \label{fig:violation_rate}
\end{figure}

\subsection{Numerical Simulations}
We tested the proposed method for the system \eqref{eq:case_dyn} with various severity of the existing failure, $\delta \in [0,10000]$, and input time delays $\tau = 0,...,4$. Our simulation results show that the generated control input $\bz^*$ is able to push the team of vehicles to form a platoon and not worsen the existing soft failure $h_0(\bx_0) \in \W_\delta$ in the future. Fig. \ref{fig:position_tau} presents examples of platoon trajectories using the control input $\bz^*$ generated by \eqref{eq:case_mpc} with various time-delay $\tau$. Since there is a non-trivial chance that the optimal solution of DRMPC may not be feasible to the original MPC problem, violations of the $\AVAR$ safety constraints will occur according to \eqref{eq:dr_prob_condition}. This phenomenon is observed in Fig. \ref{fig:violation_rate}, in which we fix the number of collected samples $N$ and solve the problem for various radii of the ambiguity set $r$ and evaluate the violation rate of safety constraints. It is intuitive to expect that the violation rate will decrease as the radius of the ambiguity set increases since it includes more potential probability distributions.

\subsection{Severity of the Existing Soft Failure}
\begin{figure}[t]
    \centering
    \includegraphics[width = \linewidth]{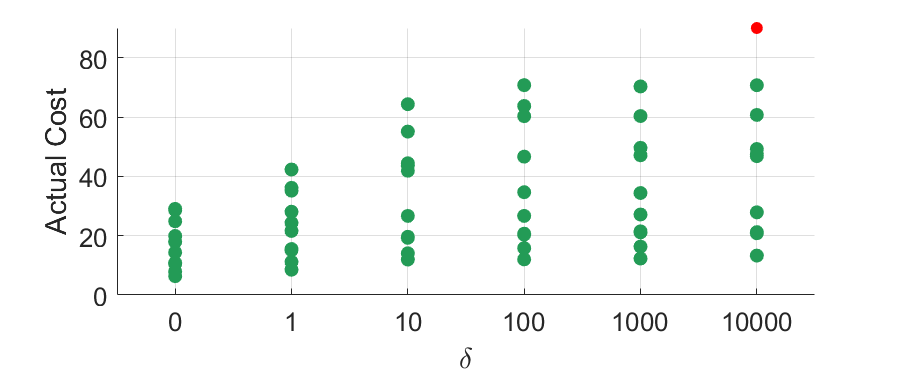}
    \caption{The actual cost the system with various severity of the existing soft failures. The red dot denotes when the problem becomes infeasible.}
    \label{fig:cost_soft}
\end{figure}

To explore how the existing soft failure $h_0(\bx_0) \in \W_\delta$ will trigger the cascading failure in the future steps, e.g., violating the safety constraints or harming the system's performance, we solve the DRMPC problem \eqref{eq:case_mpc} with different levels of $\delta$, which represents different initial inter-vehicle distances of the first pair of vehicles, i.e., $h_0(\bx_0) = x_0^{(1)} - x_0^{(2)}  = -\frac{1}{\delta+1.25}$. We evaluate the actual costs of the system \eqref{eq:case_dyn} using the optimal control input $\bz^*$ among $10$ random initial conditions and a time-delay $\tau =2$. As shown in Fig. \ref{fig:cost_soft}, for each unique $\delta$, the actual cost of the system will grow with the increase of $\delta$ since it takes more input energy and it is harder for the system to satisfy the safety constraints while minimizing the objective function. For instance, it takes more energy for the system to control a pair of vehicles that are initially close to each other to the target distance than a pair that starts at the target distance $d$ initially. Moreover, once the severity of the soft failure increases to some particular level, the problem may become infeasible since there does not exist a feasible sequence of input $\bz$ that will satisfy the $\AVAR$ safety constraints in the DRMPC problem.

\subsection{Impact of the Input Time-delay}
\begin{figure}[t]
    \centering
    \includegraphics[width = \linewidth]{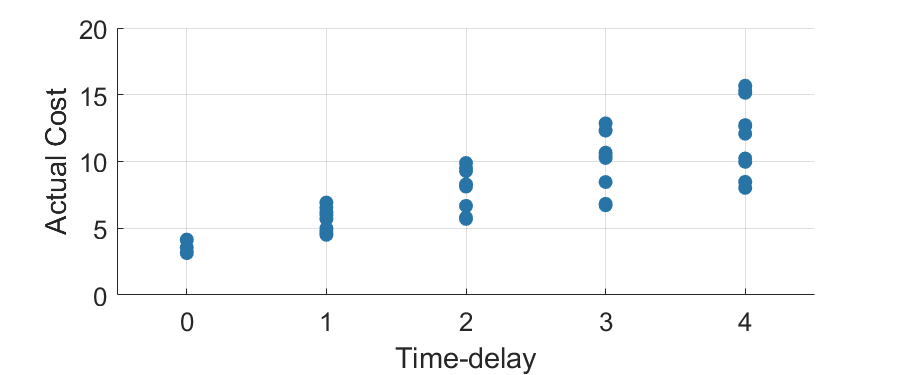}
    \caption{The actual cost the system with various input time-delays.}
    \label{fig:cost_delay}
\end{figure}

\begin{figure}[t]
    \centering
    \includegraphics[width = \linewidth]{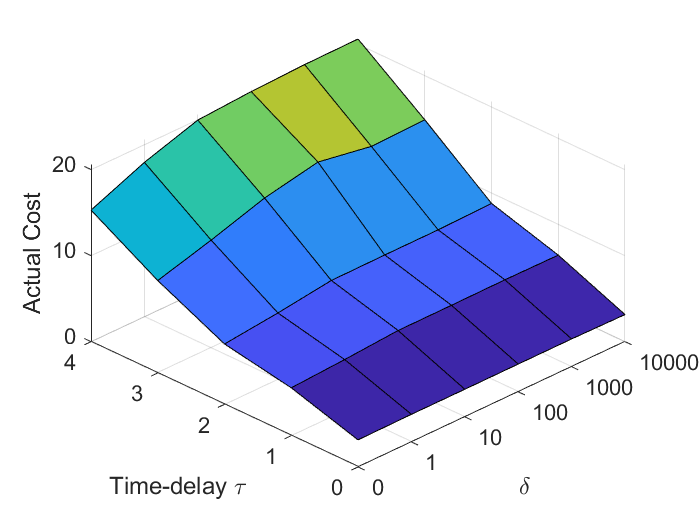}
    \caption{The actual cost of the system with various severity of the existing soft failure and time-delays.}
    \label{fig:trade_off}
\end{figure}

From equations \eqref{eq:dyn} and \eqref{eq:compact-dyn}, it becomes evident that the presence of non-zero time-delay constrains the number of effective inputs that can be applied to the system, leading to increased costs and a potential inability to achieve a feasible solution to the distributionally robust optimization problem. To demonstrate this effect, we solve the DRMPC problem \eqref{eq:case_mpc} with varying input time-delays and present our results in Fig. \ref{fig:cost_delay}. For each unique time-delay $\tau = 0,...,T-1$, we executed the optimization problem $10$ times with randomized initial conditions and applied the resulting optimal control input to the system to obtain the actual cost $J$. Our results show that the actual cost of the system increases with the increase of time-delay $\tau$, irrespective of the randomized initial conditions since the system is less likely to form a platoon from non-perfect initial conditions with the existence of the input time-delay. The collective effect of both existing soft failure and the input time-delay is shown in Fig. \ref{fig:trade_off} and it is also worth noticing that the monotonicity w.r.t severity of the existing soft failure $\delta$ and the time-delay $\tau$ in the actual cost is not guaranteed since the DRMPC program is optimizing the worst-case estimated cost $\tilde{J}$ instead of the actual cost $J$ of the system. 

\section{Conclusion}
We introduce a novel data-driven distributionally robust approach to mitigate the risk of cascading failures in delayed discrete-time LTI systems. Our proposed methodology involves formulating a distributionally robust finite-horizon optimal control problem to minimize a piecewise affine cost function while mitigating the risk of cascading failures in the presence of input time-delay and initial soft failure. Despite the challenging nature of the optimal design problem, our approach can effectively achieve the desired system performance while enforcing a set of constraints that ensure the $\AVAR$ of failures remains manageable over a finite time horizon. Our findings have significant implications for the design and implementation of control systems that are susceptible to cascading failures.

\appendix

\subsubsection{Proof of Lemma \ref{lem:1}}

When Assumption \ref{asp:light-tail} holds, using the result from measure concentration (\cite{fournier2015rate} Theorem 2), one can show that with $N$ number of samples, the Wasserstein distance between $\Pro_{\delta,\bz}$ and $\tilde{\Pro}_{\delta,\bz}$ satisfies following conditions
\begin{equation*}
    \scalebox{0.95}{$
    \Pro\left\{ W_1(\tilde{\Pro}_{\delta,\bz}, \Pro_{\delta,\bz}) \geq r \right\} \leq \begin{cases}
        c_1\exp(-c_2N r^{nT}) &\text{ if } r \leq 1\\
        c_1\exp(-c_2N r^\beta) &\text{ if } r > 1
    \end{cases},$}
\end{equation*}
where $c_1$ and $c_2$ are positive constants that only depend on $\beta, H$ and $nT$, we refer to \cite{mohajerin2018data} for detail. Then, one can conclude both branches of the result by letting $ \varepsilon = 1 - c_1\exp(-c_2N r^{nT})$ and $\varepsilon = 1 - c_1\exp(-c_2N r^\beta)$. \hfill$\square$

\subsubsection{Proof of Theorem \ref{thm:2}}

Using the result from \cite{mohajerin2018data} and \cite{micheli2022data}, and considering the fact that $c(\by,\bz)$ is a piecewise affine function, we have the worst case cost for the DRMPC problem \eqref{eq:dr_mpc},  
$
    \sup_{\Q \in \mathfrak{M}(\tilde{\Pro}_{|\delta,\bz})} \E^{\by \sim \Q} \left[ c(\by,\bz)\right]
$ 
can be obtained by solving the following convex optimization problem 
\begin{equation}    \label{eq:convex_cost}
\begin{aligned}
    \inf_{\lambda,\nu_i} ~ &\lambda \, r + \frac{1}{N} \sum_{i = 1}^{N} \nu_i\\
    \textup{subject to} ~ & a_{i_c}(\bar{A} \bz + \tilde{\bxi}^i) + b_{i_c} \bz + c_{i_c} \leq \nu_i \\
    & \|a_{i_c}\| \leq \lambda \\
    & \forall ~ i = 1,...,N \textup{ and } \forall ~ {i_c} = 1,...,N_{i_c}.
\end{aligned}
\end{equation}
Then, using the fact that $g(\by,\bz)$ is a piecewise affine function, one can convert the distributionally robust chance constraint $\sup_{\Q \in \mathfrak{M}(\tilde{\Pro}_{|\delta,\bz})} \AVAR^{\by \sim \Q}_{\alpha} \left( g(\by,\bz) \right) \leq 0$ into a convex constraint $\bz \in \mathcal{Z}_{g}$ \cite{hota2019data} , such that 
\begin{equation*}
    \scalebox{1}{$
    \mathcal{Z}_{g} := \left\{ \bz ~\Vast| ~
    \begin{aligned}
        &\theta \, r + \frac{1}{N}\sum_{i = 1}^N q_i^{(0)} \leq s^{(0)} \,\varepsilon\\
        &\max\left\{d_{i_g}(\bar{A} \bz + \tilde{\bxi}^i) + e_{i_g} + f_{i_g} + s^{(0)},0\right\}\\
        & \hspace{4cm} \leq q_i^{(0)}\\
        & \|d_{i_g}\| \leq \theta
    \end{aligned}\right\}$}
\end{equation*}
for some $s^{(0)},q_1^{(0)},...,q_N^{(0)}$, all $i = 1,...,N$ and $i_{g} = 1,..., N_{i_g}$. 

\begin{figure}[t]
    \centering
    \includegraphics[width = \linewidth]{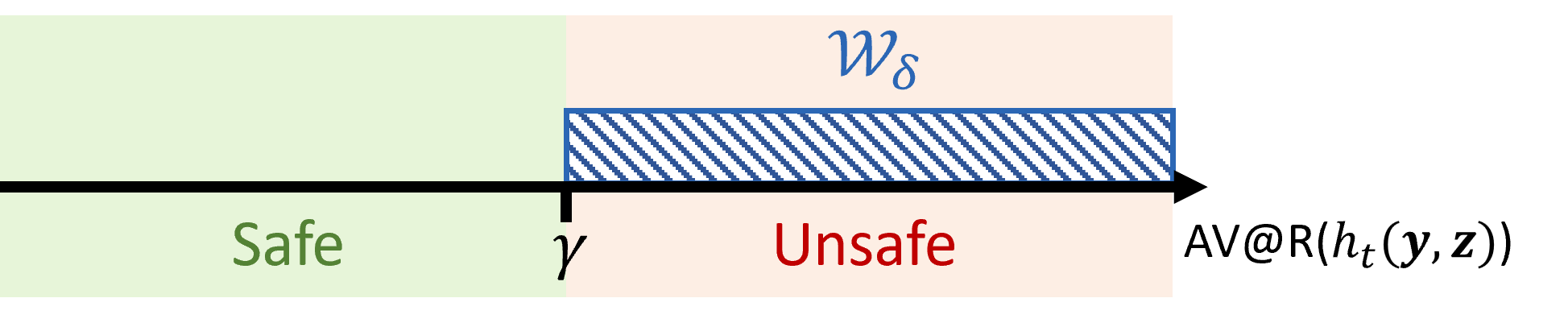}
    \caption{The diagram depicts an example of $\gamma$ for $\W_\delta \subset \R$ and the region when the $\AVAR$ constraint for $h_t(\by,\bz)$ will be satisfied (safe) or not (unsafe).}
    \label{fig:gamma}
\end{figure}

For the systemic set $\W_\delta \subset \R$, the $\AVAR$ constraints for $h_1(\by,\bz),...,h_T(\by,\bz)$ can be converted as
\begin{align*}
    \scalebox{0.95}{$
    \AVAR^{\by \sim \Q}_{\alpha} \left( h_t(\by,\bz) \right) \notin \W_\delta \Leftrightarrow \AVAR^{\by \sim \Q}_{\alpha} \left( h_t(\by,\bz) \right) - \gamma \leq 0,$}
\end{align*}
where 
\begin{equation} \label{eq:gamma}
    \gamma := ~ \inf_{h} \{h \in \R | h \in \W_\delta\}.
\end{equation}
An example of this conversion is illustrated in Fig. \ref{fig:gamma}.
Knowing that $h_t(\by,\bz)$ is a piecewise affine function $\forall t \in [0,T]$, we convert the distributionally robust chance constraint $\sup_{\Q \in \mathfrak{M}(\tilde{\Pro}_{|\delta,\bz})} \AVAR^{\by \sim \Q}_{\alpha} \left( h_t(\by,\bz) \right) \leq \gamma$ into convex constraint $\bz \in \mathcal{Z}_{h_t}$ with
\begin{equation*}
    \scalebox{1}{$
    \mathcal{Z}_{h_t} := \left\{ \bz \, \Vast| \,
    \begin{aligned}
        &\eta^{(t)}\, r + \frac{1}{N}\sum_{i = 1}^N q_i^{(t)} \leq s^{(t)} \,\varepsilon\\
        &\max\left\{\upsilon^{(t)}_{i_h}(\bar{A} \bz + \tilde{\bxi}^i) + \phi^{(t)}_{i_h} + \chi^{(t)}_{i_h} + s^{(t)},0\right\}\\
        & \hspace{4cm} \leq q_i^{(t)}\\
        & \|\phi^{(t)}_{i_h}\| \leq \eta^{(t)}
    \end{aligned}\right\}$}
\end{equation*}
for some $s^{(t)},q_1^{(t)},...,q_N^{(t)}$, all $i = 1,...,N$ and $i_{h} = 1,..., N_{i_h}^{(t)}$. Then, the DRMPC problem \eqref{eq:dr_mpc} obtains its exact convex formulation by introducing convex constraints $\bz \in \bz \in \mathcal{Z}_{g} \bigcap \mathcal{Z}_{h_1} \bigcap \cdots \bigcap \mathcal{Z}_{h_T}$ into the convex problem \eqref{eq:convex_cost}. \hfill$\square$

\printbibliography
\end{document}